\documentclass[journal ]{new-aiaa}
\usepackage[utf8]{inputenc}
\usepackage{textcomp}

\usepackage{graphicx}
\usepackage{amsmath}
\usepackage[version=4]{mhchem}
\usepackage{siunitx}
\usepackage{longtable,tabularx}

\usepackage{mathtools}
\usepackage{algorithm}
\usepackage{algpseudocode}
\usepackage{multirow}

\newcommand{\matrixstyle}[1]{\mathrm{#1}}
\newcommand{\vectorstyle}[1]{\boldsymbol{\mathbf{#1}}}

\title{Optimal Path Planning of Airborne Wind Energy Systems with a Flexible Tether}

\author{Omid Heydarnia \footnote{Doctoral Researcher, Dynamic Design Lab; email: omid.heydarnia@ugent.be}} 
\affil{Ghent University, Ghent, 9052, Belgium \\ Flanders Make, Heverlee, 3001, Belgium}
\author{Jolan Wauters\footnote{Assistant Professor; Mechatronics Group}}
\affil{KU Leuven, Bruges, 8200, Belgium \\ Flanders Make, Heverlee, 3001, Belgium}
\author{Tom Lefebvre \footnote{Assistant Professor; Dynamical Systems and Control}}
\affil{Ghent University, Ghent, 9052, Belgium \\ Flanders Make, Heverlee, 3001, Belgium}
\author{Guillaume Crevecoeur\footnote{Professor; Dynamic Design Lab}}
\affil{Ghent University, Ghent, 9052, Belgium \\ Flanders Make, Heverlee, 3001, Belgium}

\begin{document}

\maketitle

\section{Introduction}
\lettrine{T}{he} pursuit of renewable energy has driven significant research into \textsl{airborne wind energy systems} (AWES), which aim to generate electricity from more consistent and faster winds at higher altitudes. They do this by replacing the blades and towers of conventional wind turbines with an aircraft and a tether \cite{RN3}. Different types of AWES have been developed over the past decades that can be categorized based on the aircraft system (soft or rigid), generator system (fixed, moving or aircraft-mounted) and flight operation (crosswind, tether-aligned or rotational) \cite{RN26}. Particular attention has been given to a single kite (with either soft or rigid wings) that flies crosswind, unwinding a tether wound around a winch connected to a generator. \cite{RN26}. In this type of AWES, the electricity is generated through cyclic maneuvers of the aircraft, known as pumping cycles and referred to as the traction phase, during which the tether is reeled out, turning the winch and rotor of the generator. When the tether reaches its maximum length, the generator operates in the motor phase, consuming some of the harvested energy to pull the aircraft in a low-lift maneuver back to the initial point, where another pumping cycle starts. This phase is also known as the retraction phase \cite{RN3}.
        
Establishing an optimal path for the AWES to generate maximum power while keeping the aircraft safely in the sky during traction, retraction, and the transitions between these phases poses significant challenges. This can be attributed to the strong non-linearity of the dynamics and underlying aerodynamics, and the presence of the tether. A vast body of work has addressed these challenges in the pursuit of feasibility, accuracy, and computational efficiency. Gros \& Diehl proposed a non-minimal coordinate representation model considering the tether as a rigid rod and modeling the aircraft and tether with an index-3 \textsl{differential-algebraic system of equations} (DAE) \cite{RN30}. This approach was also followed by \cite{RN15,RN13,RN5}. On the other hand, Erhard et al. utilized a quaternion-based representation in an \textsl{ordinary differential system of equations} (ODE), assuming the tether as a rigid rod, with its length varying with winch velocity and its force being estimated based on airspeed and aerodynamic characteristics of the kite \cite{RN31}. Simpler modeling approaches also assuming a rigid tether, were pursued by \cite{ilzhofer2007nonlinear,RN14,RN12}. 
       
Although it is a valid assumption to consider the tether as a rigid rod during the traction phase as it is under high tension \cite{RN8}, this assumption is harder to make in the retraction and transition phases, when sag of tether might occur. 


    Various studies have been conducted to develop more accurate models for tethers. Du et al. derived the equations for a variable-length cable using the absolute nodal coordinate formulation and discretized the cable into a finite number of elements to reduce computational cost \cite{du2015dynamic}. Similarly, \cite{RN28} utilizes the same approach to model the cable, and solves the whole AWES, including aircraft, tether and winch system, by a variational integrator to increase the performance. In another approach, the tether is considered to correspond to a fixed number of lumped masses connected by spring-damper elements to generate a coherent model \cite{RN27,berra2019optimal}. Although modeling the tether with the above-mentioned methods increases the accuracy, they add three degrees of freedom for each rod/lumped mass, increasing the whole system's number of dimensions drastically. To tackle this difficulty, \cite{RN19} introduced a quasi-static approach, in which the elastic vibrations are neglected and the cable tension and shape are assumed to be the results of mass and drag. Although the simulation of small-scale AWES using a quasi-static approach for the tether has been verified through flight tests \cite{williams2019flight}, only one study exists that uses a flexible tether model for path planning \cite{RN36}. In that study only the landing of a simple aircraft in 2D is examined.

    The major contribution of this note is considering the flexible tether model in an optimal control problem such that the computational cost remains manageable. Then, the obtained results of the flexible tether are contrasted with the rigid tether model in terms of power and tether force magnitude and direction at the aircraft. From a numerical optimization standpoint, we introduce an index-1 DAE formulation for AWES with a minimal coordinate representation, which eliminates consistency condition problems and drift arising from the index reduction of higher index DAE formulations. The problem is transcribed into a direct multiple shooting formulation using the penalty-based homotopy approach, whose efficiency in increasing the convergence rate and improving the solver stability has been proven \cite{RN5}. 
    
    The note is structured as follows, in section \ref{sec_2} we describe the AWES model which builds on the MegAWES toolbox \cite{RN17}. Subsequently, the optimization method and its implementation using CasADi \cite{RN68} are discussed in section \ref{sec_3}. Finally, in section \ref{sec_4}, the simulation results are presented to show the effectiveness of the proposed optimization method for AWES with the flexible tether.

\section{Modeling}
\label{sec_2}
    
    Our modeling approach follows the methodology established by MegAWES \cite{RN17}. The model distinguishes itself from those used by other optimization tools, such as AWEBOX \cite{RN5}, in two critical aspects.
    First, a minimal coordinate representation is used to parameterize the aircraft's attitude instead of an element of $\mathfrak{SO}(3)$. This is strictly a modeling choice. Second, a flexible tether model is included instead of a rigid tether model. A flexible tether increases the physical accuracy of the resulting model. Both of these model features have a positive effect on the subsequent optimization. This will be discussed in more detail in the next section.

    \subsection{Parameterization and kinematics}
    Consider two frames of reference $\mathcal{O}$ and $\mathcal{B}$. Reference frame $\mathcal{O}$ is used as an inertial reference frame, whereas frame $\mathcal{B}$ is attached to the aircraft. Therefore, the attitude of the aircraft can be represented as the rotation matrix $\matrixstyle{R}_{\mathcal{O}}^\mathcal{B}\in\mathfrak{SO}(3)$, that maps $\mathcal{O}$ into $\mathcal{B}$. As is standard in the aerospace community, we use the Euler angles to parameterize the attitude\footnote{For optimization purposes, other conventions can be considered to avoid singularities.}. The angles are gathered in a coordinate vector, $\vectorstyle{q}_\text{a}\in\mathbb{R}^3$.
	\begin{equation}
		\matrixstyle{R}_{\mathcal{O}}^\mathcal{B} = \matrixstyle{R}_x(\phi) \matrixstyle{R}_y(\theta) \matrixstyle{R}_z(\psi)
 	\end{equation}
 	
 	The time derivative of the coordinate vector $\vectorstyle{q}_\text{a}$ is related to the angular velocity of the aircraft by a geometric Jacobian 
 	\begin{equation}
            \label{eq:euler}
 		\vectorstyle{\omega}_{\text{aircraft}}^\mathcal{B} = \begin{pmatrix}
 			c_\theta c_\psi & s_\psi & 0 \\
 			-s_\psi & c_\psi & 0 \\
 		s_\theta c_\psi & 0 & 1
 		\end{pmatrix}\dot{\vectorstyle{q}}_\text{a}
 	\end{equation}
   where $c_{(.)}$ and $s_{(.)}$ denote $\cos{(.)}$ and $\sin{(.)}$, respectively. 
   
   We introduce a third frame of reference $\mathcal{W}$, whose $x$-axis aligns with the wind direction at an angle $\zeta$ from the inertial $x$-axis, and whose $z$-axis points upward, opposite to the downward direction of the inertial $z$-axis.
	\begin{equation}
		\matrixstyle{R}_{\mathcal{O}}^\mathcal{W} =  \matrixstyle{R}_z(\zeta) \matrixstyle{R}_x(\pi)
	\end{equation}
	
\begin{figure*}[ht!]
        \centering
        \includegraphics[trim={2.0cm 1.5cm 2.0cm 1.05cm}, clip, width=.6\textwidth]{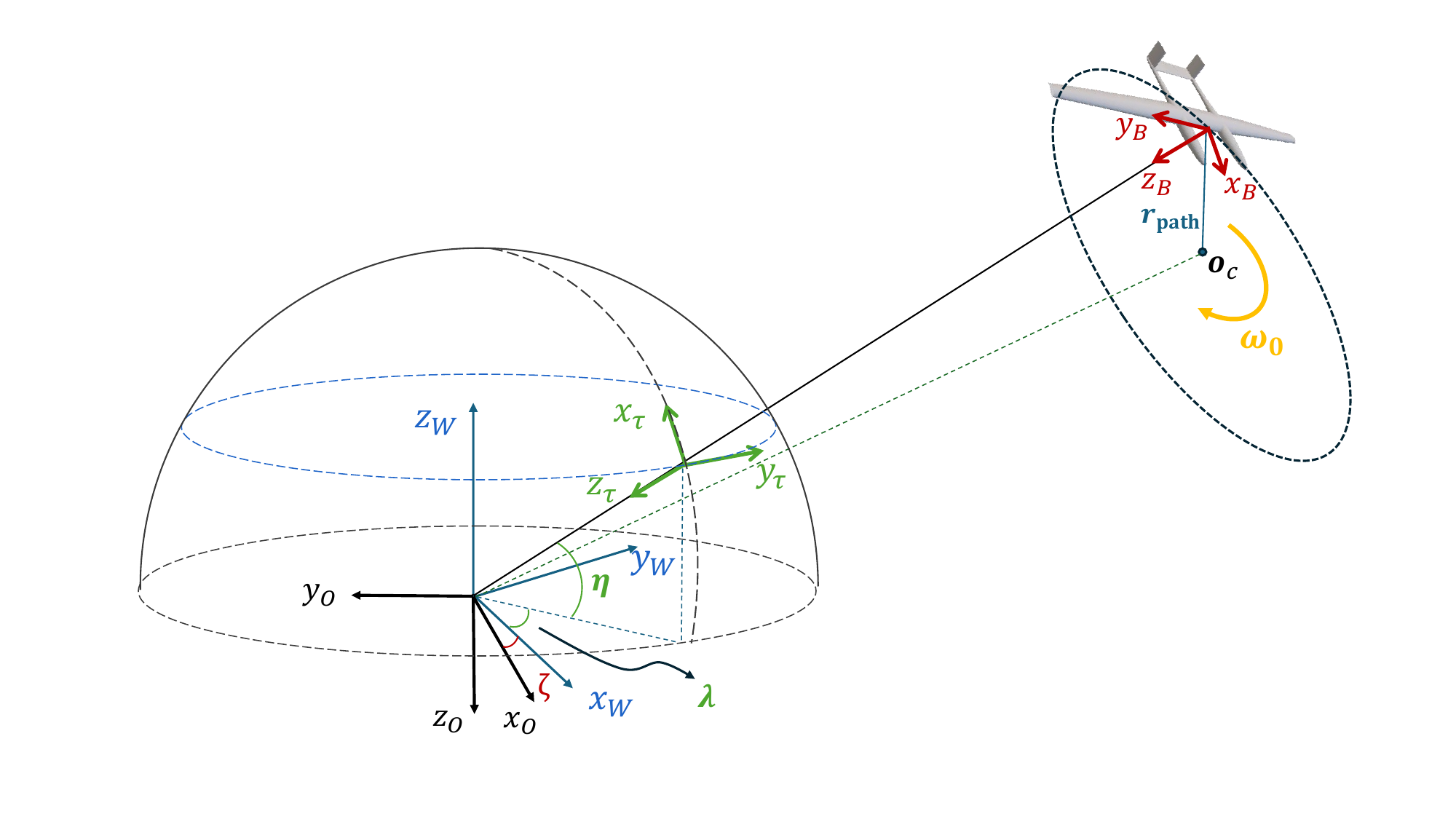}
        \caption{Visualization of initial circular path concept, and the reference frames ($\mathcal{O}$ (black), $\mathcal{W}$ (blue), $\mathcal{\tau}$ (green), $\mathcal{B}$ (red)).}
        \label{fig:path_initialization}
    \end{figure*}

	The position of the aircraft is parameterized using spherical coordinates $\vectorstyle{q}_\text{s}\in\mathbb{R}^3$, where the first and second elements of $\vectorstyle{q}_\text{s}$ are the azimuth and polar angles, as shown in Fig.~\ref{fig:path_initialization}, and the third one denotes the radial distance from the origin to the aircraft. To that end, we introduce a forth and final frame of reference $\mathcal{T}$, whose motion relative to $\mathcal{W}$ is parameterized by the azimuth $\lambda$, and polar $\eta$ angles and so that the local $z$-axis is pointing away from the aircraft to the origin. In local coordinates, the position of the aircraft is then characterized by the distance $r$ along the $z$-direction.
	\begin{equation}
		\matrixstyle{R}_{\mathcal{W}}^\mathcal{T} = \matrixstyle{R}_z(\lambda)\matrixstyle{R}_y(-\eta) \matrixstyle{R}_y(-\tfrac{\pi}{2})
	\end{equation}
	
	The time derivative of the coordinate vector $\vectorstyle{q}_\text{s}$ is related to the linear velocity of the aircraft by another geometric Jacobian.
	\begin{equation}
            \label{eq:position}
		      \vectorstyle{v}_{\text{aircraft}}^\mathcal{W} = \begin{pmatrix}
		 	c_\lambda c_\eta & -s_\lambda c_\eta r & -c_\lambda s_\eta r \\
		 	s_\lambda c_\eta & c_\lambda c_\eta r & -s_\lambda s_\eta r \\
		 	s_\eta & 0 & c_\eta r 
		 \end{pmatrix}\dot{\vectorstyle{q}}_\text{s}
	\end{equation}
    A visualization of the frames is presented in Fig.~\ref{fig:path_initialization}.
   
    \subsection{Rigid-body aircraft dynamics}
	The rigid body dynamics of an aircraft are described in numerous references, such as \cite{doi:https://doi.org/10.1002/9781119174882.ch2}. The main difference between modeling an ordinary aircraft and a tethered one, is the tether force. In this context, we assume that the tether force acts on the aircraft's center of gravity. This means that no torque is applied to the aircraft due to the tether. The aircraft's dynamics are governed by
        \begin{subequations}
	   \label{eq:f1}
	   \begin{align}
		\label{eq:f1a}
		      m  \dot{\vectorstyle{v}}_{\text{aircraft}}^\mathcal{B} + m \vectorstyle{\omega}_{\text{aircraft}}^\mathcal{B}\times \vectorstyle{v}_{\text{aircraft}}^\mathcal{B} &= \vectorstyle{f}_{\text{aircraft}}^\mathcal{B} \\
     \label{eq:f1b}
		\matrixstyle{J}\dot{\vectorstyle{\omega}}_{\text{aircraft}}^\mathcal{B} +  \vectorstyle{\omega}_{\text{aircraft}}^\mathcal{B}\times \matrixstyle{J}\vectorstyle{\omega}_{\text{aircraft}}^\mathcal{B}&= \vectorstyle{\tau}_{\text{aircraft}}^\mathcal{B}
	\end{align}
    \end{subequations}
    where $\vectorstyle{v}_{\text{aircraft}}^\mathcal{B}$ is the aircraft velocity in the body reference frame, $\vectorstyle{\omega}_{\text{aircraft}}^\mathcal{B}$ is the angular velocity vector between the body and inertial reference frames. $m$ and $\matrixstyle{J}$ represent the aircraft mass and inertia, respectively.
    
    The total force, $\vectorstyle{f}_{\text{aircraft}}$, includes aerodynamic, tether, and gravitational forces. The total torque, $\vectorstyle{\tau}_{\text{aircraft}}$, includes only aerodynamics torques. 
    \begin{subequations}
        \label{eq:aircraft}
        \begin{align}
            \label{eq:forces}
            \vectorstyle{f}_{\text{aircraft}}^\mathcal{B} &= \vectorstyle{f}_{\text{aircraft,aero}}^\mathcal{B} + \vectorstyle{f}_{\text{aircraft,tether}}^\mathcal{B} + 
            \matrixstyle{R}_{\mathcal{O}}^\mathcal{B} (mg \vectorstyle{E}_3^\mathcal{O})  \\
            \label{eq:torques}
            \vectorstyle{\tau}_{\text{aircraft}}^\mathcal{B} &=       \vectorstyle{\tau}_{\text{aircraft,aero}}^\mathcal{B}
            \end{align}
	\end{subequations}

    \subsection{Aerodynamic model}
    In this work, we adopt the MegAWES aerodynamic model \cite{RN17}. MegAWES utilizes pre-calculated static aerodynamic coefficients stored in lookup tables to compute the aerodynamic forces and torques acting on the wing, elevator, and rudder, independently, as a function of the aircraft's apparent velocity $\vectorstyle{v}_a$ (\ref{eq:V_a}), angle of attack $\alpha$ (\ref{eq:AoA}), side slip angle $\beta$ (\ref{eq:sl}), and deflections of control surfaces $\delta$.
    \begin{subequations}
	\begin{align}
        \label{eq:V_a}
        \vectorstyle{v}_a &= \vectorstyle{v}_\text{aircraft}^\mathcal{W} - \vectorstyle{v}_\text{wind}^\mathcal{W} \\
        \label{eq:AoA}
        \alpha &= -\arctan \frac{v_\text{aircraft,z}^\mathcal{B}-v_\text{wind,z}^\mathcal{B}}{v_\text{aircraft,x}^\mathcal{B}-v_\text{wind,x}^\mathcal{B}}\\
        \label{eq:sl}
		\beta &= \arctan \frac{v_\text{aircraft,y}^\mathcal{B}-v_\text{wind,y}^\mathcal{B}}      {v_\text{aircraft,x}^\mathcal{B}-v_\text{wind,x}^\mathcal{B}}
	\end{align}
    \end{subequations}
   
    To retain a differentiable model, the look-up tables have been replaced with second-order polynomials of the following form
        \begin{equation}
            \begin{aligned}
                \label{eq:aero_coefs}
                    \mathcal{C}_{i}(\alpha) &= \mathcal{C}_{i,2} \alpha^2 + \mathcal{C}_{i,1} \alpha + \mathcal{C}_{i,0}     
            \end{aligned}
        \end{equation}
    where the index $i \in\{L, D, M\}$ determines the coefficients of lift and drag forces, and the pitch moment of the wing. Similarly, for the rudder and elevator, $i \in\{L, D\}$ specifies the lift and drag forces. Additional information can be found in \cite{RN17}.

     \subsection{Winch model}
	AWES with a fixed ground station consist of a winch and a generator. We consider a simple 1-DOF winch model. This is similar to \cite{RN18}.
	\begin{equation}
			\label{eq:f2}
		J_{\text{winch}}\ddot{\theta}_{\text{winch}} + b_{\text{winch}} \dot{\theta}_{\text{winch}} = r_\text{winch} \|\vectorstyle{f}_{\text{winch,tether}}\| - \tau_{\text{generator}}
	\end{equation}
    where $\dot{\theta}_{\text{winch}}$, $J_{\text{winch}}$, $r_\text{winch}$, and $b_{\text{winch}}$ represent the winch's angular velocity, inertia, radius, and viscous friction, respectively. The amplitude, $\|\vectorstyle{f}_{\text{winch,tether}}\|\neq0$, is equal to the tether force magnitude evaluated at the winch. Finally we have the the motor/generator torque, $\tau_{\text{generator}}$, which will be treated as a control input.

    \subsection{Quasi-static flexible tether model}

    Other references \cite{RN13,RN5,RN31} simplify the tether physics to ease optimization, but this sacrifices accuracy, which can significantly affect the results. This work aims to integrate a tether model that balances physical accuracy with computational efficiency. To achieve this, we adopt the quasi-static model from  \cite{RN19}, which accounts for tether sag and elastic deformation. While it does not capture dynamic effects, it provides a physically accurate and efficient approximation using discrete lumped masses connected by massless elastic (or inelastic) segments. 

    The tether length is a dynamic variable that is related to the other system variables.
    \begin{equation}
		\label{eq:theta_winch}
		l_\text{tether} =  r_\text{winch} \theta_\text{winch}
    \end{equation}

    We approximate the cable with a total of $N+1$ lumped masses. We describe every mass, $m_j$, with its position which is denoted with the vector $\vectorstyle{p}_j$. Indexing is done such that $\vectorstyle{p}_0$ coincides with the aircraft and $\vectorstyle{p}_{N}$ with the winch. The nominal length and mass of each cable segment are defined as
	\begin{subequations}
		\begin{align}
			L_s &= \tfrac{1}{N+1} l_\text{tether}\\
			m_j &=  \rho_{\text{tether}} L_s
		\end{align}
	\end{subequations}
    where $\rho_{\text{tether}}$ represents the cable's mass density.
 
    The principle idea behind the model is to construct a recursive calculation that iterates over all masses in the discretized cable model and determines their positions and associated constraint forces. The model takes the tether force at the lumped mass $f_N$ on the winch as input and outputs the corresponding position of the final mass $p_0$, which should coincide with the aircraft's position. Once the aircraft's position is known, a root-finding problem can be formulated to determine the corresponding tether force.
    
    To determine the calculation procedure, we can consider the equation of motion of the $j$\textsuperscript{th} mass. To simplify the notation the superscript is left out in this section and all vectors are expressed in the wind frame of reference
	\begin{equation}
            \label{eq:tether_forces}
		m_j \ddot{\vectorstyle{p}}_j = \vectorstyle{f}_{j-1} - \vectorstyle{f}_{j} + \vectorstyle{f}_{j,\text{drag}}- m_j g \vectorstyle{E}_3
	\end{equation}
    where $\vectorstyle{f}_j$ denotes the constraint force shared with the preceding mass,  $\vectorstyle{f}_{j-1}$ denotes the constraint force shared with the following mass,  $\vectorstyle{f}_{j,\text{drag}}$ denotes the local drag force and the final term represents gravitational effects. We can now further specify all elements from the equation. The recursion will follow.
    
    First it is assumed that the velocity of each mass along the cable is computed as the sum of the radial velocity induced by the aircraft and the angular rotation of the segment and the the acceleration of each segment is computed from the centrifugal component only. In other words, the tether behaves as a rigid body at every time instant.
	\begin{subequations}
            \begin{align}
	           \dot{\vectorstyle{p}}_j &= \vectorstyle{v}_\text{tether} + \vectorstyle{\omega}_{\text{tether}}  \times \vectorstyle{p}_j  \\
			\ddot{\vectorstyle{p}}_j &= \vectorstyle{\omega}_{\text{tether}}  \times 	           (\vectorstyle{\omega}_{\text{tether}}  \times \vectorstyle{p}_j)
            \end{align}
	\end{subequations}

    The radial velocity and angular rotation of the aircraft are defined as
    \begin{subequations}
		\begin{align} 
			\vectorstyle{v}_{\text{tether}} &=  \tfrac{1}{\| \vectorstyle{p}_{\text{aircraft}}
				\|^2} \left\langle \vectorstyle{p}_{\text{aircraft}},\vectorstyle{v}_{\text{aircraft}}\right\rangle  \vectorstyle{p}_{\text{aircraft}} \\
		\vectorstyle{\omega}_{\text{tether}} &= \tfrac{1}{\|\vectorstyle{p}_{\text{aircraft}}
			\|^2} \vectorstyle{p}_{\text{aircraft}}\times \vectorstyle{v}_{\text{aircraft}}
	\end{align}
	\end{subequations}
	
    The drag force acting on the $j$\textsuperscript{th} mass is defined as
    \begin{equation}
    	\vectorstyle{f}_{j,\text{drag}}= - \frac{1}{2} \rho d_\text{tether} l_j C_D \|\vectorstyle{v}_{j,\perp}\|\vectorstyle{v}_{j,\perp}
    \end{equation}
    where $\rho$ is the air density at the $j$\textsuperscript{th} segment, $d_\text{tether}$ represents the tether diameter, $l_j$ is the unstrained length of $j$\textsuperscript{th} segment, $C_D$ is the drag coefficient, and $\vectorstyle{v}_{j,\perp}$ is the normal component of relative wind velocity at each segment. The normal  $\vectorstyle{v}_{j,\perp}$ and parallel $\vectorstyle{v}_{j,\parallel}$  components of each segment can be calculated as
        \begin{subequations}    
            \begin{align}
                    \label{eq:v_prep}
        		\vectorstyle{v}_{j,\perp} &= \vectorstyle{v}_{j} - 	\vectorstyle{v}_{j,\parallel} \\
        		\vectorstyle{v}_{j,\parallel} &= \tfrac{1}{\|\Delta \vectorstyle{p}_{j}
                    \label{eq:v_parr}
        			\|^2} \left\langle\Delta \vectorstyle{p}_{j},\vectorstyle{v}_j\right\rangle \Delta \vectorstyle{p}_{j} 
            \end{align}
        \end{subequations}
    where $	\Delta \vectorstyle{p}_j =  \vectorstyle{p}_{j-1} -  \vectorstyle{p}_{j}$, and relative wind velocity at each segment reads as
        \begin{equation}
        	\vectorstyle{v}_j = \dot{\vectorstyle{p}}_j  - \vectorstyle{v}_{j,\text{wind}} = \dot{\vectorstyle{p}}_j - \vectorstyle{v}_{\text{wind}}(z_j)
        \end{equation}
    
    \begin{figure}[ht!]
        \centering
        \includegraphics[trim={3cm 1.0cm 3.8cm 1.2cm},clip, width=.5\textwidth]{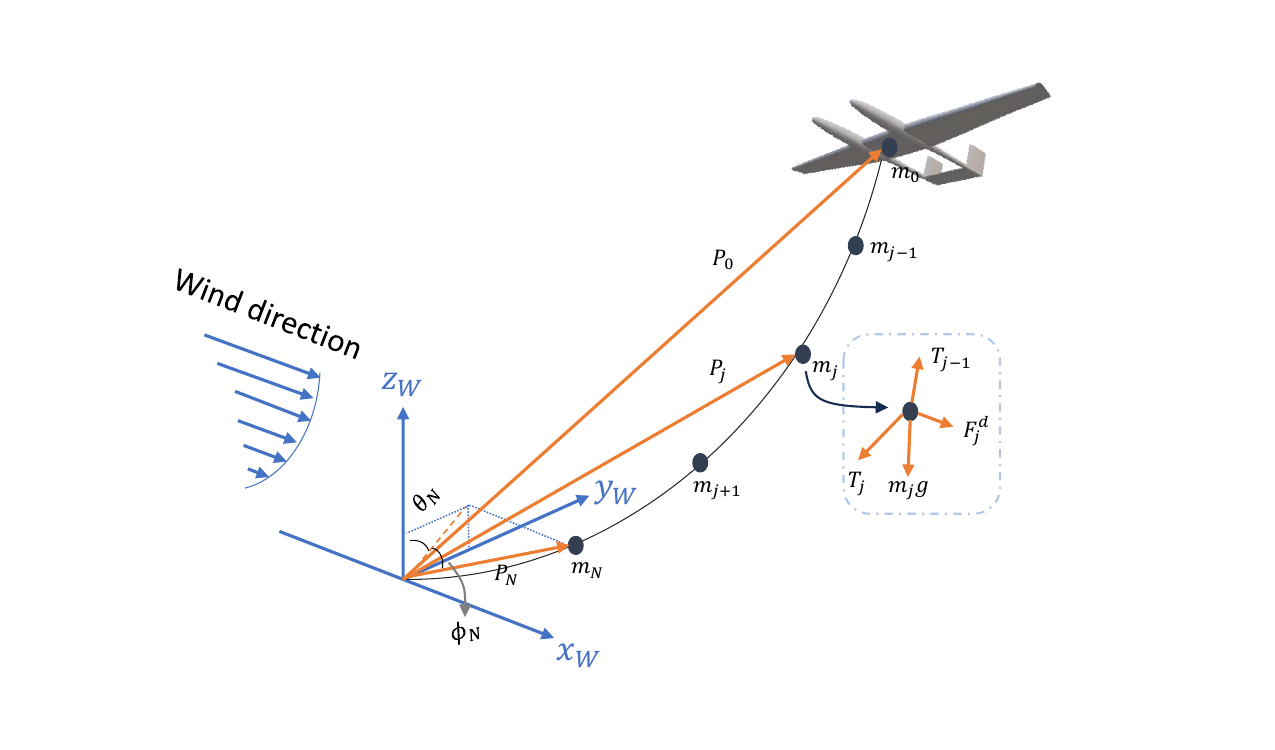}
        \caption{Visualization of tether lumped masses and forces in the wind frame.} 
        \label{fig:tether}
    \end{figure}

    Now provided that $\vectorstyle{f}_j$ and $\vectorstyle{p}_j$ are known, equation (\ref{eq:tether_forces}) can be manipulated to yield $\vectorstyle{f}_{j-1}$. Once $\vectorstyle{f}_{j-1}$ is known, the position of the $j-1$\textsuperscript{th} mass can be calculated. According to the Hook's law we have that
        \begin{equation}
            l_{j-1} = L_s \left(1+\frac{\|\vectorstyle{f}_{j-1}\|}{E_{\text{tether}} A_\text{tether}}\right)
        \end{equation}
    where $E_{\text{tether}}$ and $A_{\text{tether}}$ are tether Young's modulus and cross-section area. Due to the assumption that the tether is in static equilibrium it also follows that
        \begin{equation}
    			\vectorstyle{p}_{j-1} = \vectorstyle{p}_j + l_{j-1} \frac{\vectorstyle{f}_{j-1}}{\|\vectorstyle{f}_{j-1}\|}
        \end{equation}

    These equations thus produce the anticipated recursion. We can parameterize the first constraint force and solve for $\vectorstyle{p}_0$. By definition, we find the following algebraic equation
    \begin{equation}
        \begin{aligned}
        \label{eq:tether_algebraic}
        \vectorstyle{p}_{\text{aircraft}}^\mathcal{W} - \vectorstyle{p}_0 = 0
        \end{aligned}
   \end{equation}
   
   Solving for the roots of this equation as a function of the first constraint force determine the quasi-static model. The force acting on the aircraft is evaluated as 
   \begin{equation}
        \begin{aligned}
            \vectorstyle{f}_{\text{aircraft}, \text{tether}}^\mathcal{B} = -\matrixstyle{R}_{\mathcal{O}}^{\mathcal{B}} \matrixstyle{R}_{\mathcal{W}}^{\mathcal{O}} \vectorstyle{f}_0
        \end{aligned}
   \end{equation}

   In this work we parameterize the first constraint force as follows   
        \begin{subequations}
        \label{eq:tether_nPoint}
        \begin{align}
        \label{eq:tether_nPointa}
        \vectorstyle{f}_N &= T \begin{pmatrix} s_{\theta_N} c_{\phi_N} & s_{\phi_N} & c_{\theta_N}c_{\phi_N} \end{pmatrix}^\top \\
        \label{eq:tether_nPointb}
        \vectorstyle{p}_N &= L_s\begin{pmatrix} s_{\theta_N}c_{\phi_N} & s_{\phi_N} & c_{\theta_N}c_{\phi_N} \end{pmatrix}^\top
        \end{align}
   \end{subequations}
   where $T$, $\theta_N$ and $\phi_N$ denote the magnitude of the force and the spherical angles that determine the direction of force to the nearest point mass to the winch, as depicted in Fig.~\ref{fig:tether}. 

\section{Optimization}
    \label{sec_3}
    The goal of this work is to integrate the model from the previous section into an optimization framework tailored to generating periodic flight trajectories that maximize power generation. This means that we are looking for feasible trajectories that can be repeated infinitely and convert the maximum amount of wind energy into electrical energy with every cycle. This problem formulation results into a challenging optimal control problem that requires careful numerical treatment. Our approach draws much inspiration from AWEBOX \cite{RN5}, but instead uses a minimal coordinate representation of the attitude and a flexible tether model. The use of the quasi-static flexible tether model furthermore introduces some additional challenges that require further treatment. In this section, we first define all the components of the optimization problem. We then describe our solution strategy, which has been specifically tailored to meet the unique requirements of this problem.
    
    \subsection{Problem formulation}
    We aim to generate periodic maximum power flight trajectories for AWES. Similar to other work, we propose to solve the following periodic continuous-time Optimal Control Problem (OCP). 
    \begin{equation}
            \label{eq:ocp}
            \begin{aligned}
			&\min_{\vectorstyle{w}} \frac{1}{T} \int_{t=0}^{T} c(\vectorstyle{x}     
                (t),\vectorstyle{u}(t),\vectorstyle{z}(t))\text{d}t \\
			&\begin{aligned}\text{s.t. }
					 0&= \vectorstyle{F}(\dot{\vectorstyle{x}}(t),\vectorstyle{x}(t),\vectorstyle{u}(t),\vectorstyle{z}(t)), \forall t \in [0,T) \\
					0 &\geq  \vectorstyle{H}(\vectorstyle{x}(t),\vectorstyle{u}(t),\vectorstyle{z}(t)), \forall t \in [0,T) \\
					0&= \vectorstyle{x}(0) -\vectorstyle{x}(T) \\
					0&= \Psi(\vectorstyle{x}(0))	\end{aligned} \\
            \end{aligned}
	\end{equation}

    The OCP is fully specified by definition of the cost rate function $c$, the differential and possibly other equality constraints $\vectorstyle{F}$, and the inequality path constraints $\vectorstyle{H}$.
    
    Here, the optimization parameters are denoted by $ \vectorstyle{w} = \begin{pmatrix} \vectorstyle{x}^\top & \vectorstyle{u}^\top & \vectorstyle{z}^\top\end{pmatrix}^\top $, which represent the concatenation of the differential states $\vectorstyle{x}$, control inputs $\vectorstyle{u}$, and algebraic states $\vectorstyle{z}$, defined as 
    \begin{subequations}
		\begin{align}
            \label{eq:diffrential_states}
		\vectorstyle{x} &= \begin{pmatrix}
		\theta_{\text{winch}} &
		\dot{\theta}_{\text{winch}}&
            \vectorstyle{v}_{\text{aircraft}}^{\mathcal{B},\top} &
		\vectorstyle{\omega}_{\text{aircraft}}^{\mathcal{B},\top} &
            \vectorstyle{q}_\text{a}^\top  &
            \vectorstyle{q}_\text{s}^\top
		\end{pmatrix}^\top \\
            \label{eq:algebraic_states}
		\vectorstyle{z} &= \begin{pmatrix}
            \phi_N &
            \theta_N &
            T
            \end{pmatrix}^\top \\
            \label{eq:ctrl_inputs}
		\vectorstyle{u} &= \begin{pmatrix} 
            \tau_{\text{generator}}&
            \delta_{\text{aileron}} &                 
            \delta_{\text{elevator}} &
		\delta_{\text{rudder}}        
		\end{pmatrix}^\top 
	\end{align}
	\end{subequations}
    
    The OCP in (\ref{eq:ocp}) makes use of a general form of DAEs defined in fully implicit form, $\vectorstyle{F}$. In this work a semi-explicit DAE description of the system dynamics is available so that we can further specify the DAE as    \begin{equation}
        \label{eq:equality}
	\begin{aligned}
            \vectorstyle{F}={\begin{pmatrix}
                \vectorstyle{\dot{x}}-\vectorstyle{f}(\vectorstyle{x}(t),\vectorstyle{u}(t),\vectorstyle{z}(t)) \\
                \vectorstyle{g}(\vectorstyle{x}(t),\vectorstyle{z}(t))
            \end{pmatrix}} = 0
        \end{aligned}
    \end{equation}
    where the algebraic $\vectorstyle{g}$, and differential $\vectorstyle{f}$, equations corresponded to the tether dynamics (\ref{eq:tether_algebraic}) and the concatenation of the winch (\ref{eq:f2}) and aircraft dynamics (\ref{eq:euler}, \ref{eq:position}, \ref{eq:f1}), respectively. 

    Similar to other studies such as \cite{RN13, RN5}, the cost function in (\ref{eq:ocp}) can be expressed as the sum of negative power output and a penalty on the angular velocity rate of the aircraft and side slip angle, to maximize the power output whilst also avoiding aggressive maneuvers. 
	\begin{equation}
        \label{eq:cost}
		c(\vectorstyle{x},\vectorstyle{u},\vectorstyle{z}) = - P+ \|\vectorstyle{\xi}\|_\Gamma^2
	\end{equation}
    Here \(P= \dot{\theta}_{\text{winch}} \tau_{\text{generator}}\), \(\vectorstyle{\xi} = \begin{pmatrix} \dot{\vectorstyle{\omega}}_{\text{aircraft}}^\mathcal{B} & \beta \end{pmatrix}^\text{T}\), and $\Gamma$ is a constant diagonal weight matrix.
    
\begin{figure}[ht!]
                \centering
                \includegraphics[trim={0.5cm 1.05cm 0.5cm 1.05cm}, clip, width=.5\textwidth]{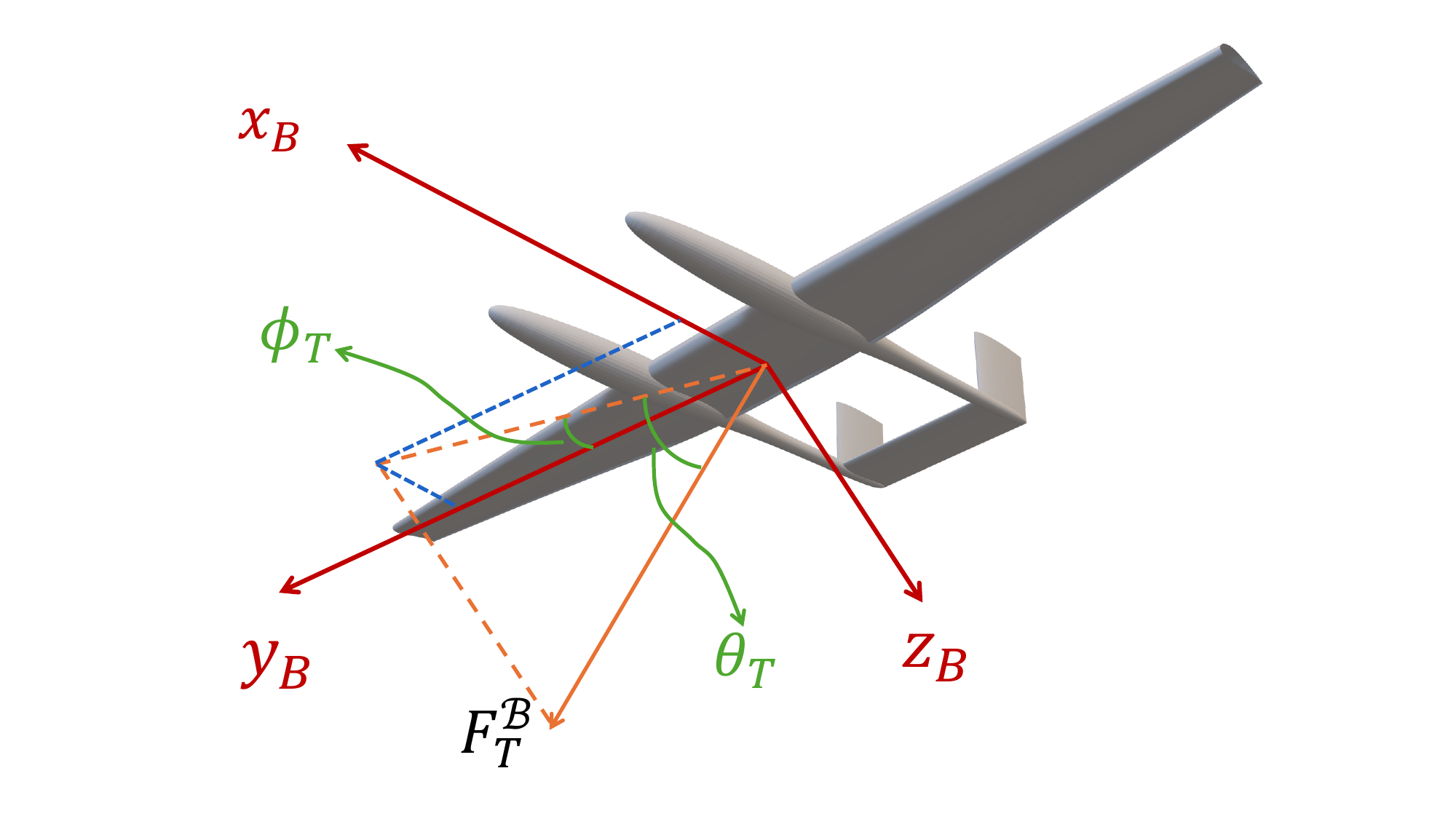}
                \caption{The constraint on tether angle to avoid collision.} 
                \label{fig:tether_collision}
            \end{figure}

    The inequality constraints or path constraints $\vectorstyle{H}$ plays a critical role in generating realistic trajectories. Constraints such as the maximum force that tether can tolerate, flight region, winch acceleration, and aerodynamics and control surfaces limitations must be taken into account for trajectory generation. Furthermore, the tether and aircraft must not collide. Accordingly, the angle between $\vectorstyle{f}_{\text{tether}}^\mathcal{B}$ and the $xy$-plane in the body reference frame $\theta_T$ is restricted to avoid any contact between tether and aircraft, where the angles describing the tether force in the body reference frame are chosen as shown in Fig.~\ref{fig:tether_collision}.    
    \begin{equation}
        \label{eq:tether_collision}
        \begin{aligned}
            \underline{\theta_T} \le \theta_T \le \overline{\theta_T}
        \end{aligned}
    \end{equation}
    
    The initial value constraint is embedded in function $\Psi$, where $\Psi(\vectorstyle{x})= 0$ is used to enforce consistency conditions for higher-order DAEs. Here, it can be utilized for the cold start of the tether dynamics and for defining other constraints that can enhance the performance of the optimal path or might assist the system operation in take-off or landing. These aspects are further addressed in the simulation and results section.
 
    So far, all elements of the OCP problem (\ref{eq:ocp}) have been introduced. Solving (\ref{eq:ocp}) is non-trivial because the DAE $\vectorstyle{F}$, defines a constraint manifold that contains all achievable trajectories of the system. For strongly under-actuated and highly non-linear dynamics, this manifold can become highly complex, making it increasingly difficult to find feasible solutions let alone navigate the manifold to identify the optimal trajectory. As a result, the optimization algorithm may converge prematurely or fail to find a solution altogether. To ease the challenges faced by the optimization algorithm, a homotopy approach is used here.
    
    \subsection{Homotopy based approach}
     The purpose of homotopy based strategies is to decrease the overall computational cost and improve the solver's reliability and robustness \cite{RN47,RN5,he2022homotopy, bergman2018combining}. Within a homotopy approach, challenging optimization problems are treated by solving a sequence of easier, manageable problems that gradually evolve to the original, complex problem. To that end, the OCP in (\ref{eq:ocp}) is reformulated as
    
    \begin{equation}
            \label{eq:ocp_homotopy}
            \begin{aligned}
			&\min_{\tilde{\vectorstyle{w}}, \vectorstyle{\Phi}} \frac{1}{T} \int_{t=0}^T c_\Phi(\vectorstyle{x}     
                (t),\tilde{\vectorstyle{u}},\vectorstyle{z}(t), \vectorstyle{\Phi}) + \vectorstyle{s}^\top\vectorstyle{\Phi} \;\text{d}t\\
			&\begin{aligned}\text{s.t. }
					 0&= \vectorstyle{F}_{\Phi}(\dot{\vectorstyle{x}}(t),\vectorstyle{x}(t),\tilde{\vectorstyle{u}}(t),\vectorstyle{z}(t),\vectorstyle{\Phi}), ~ \forall t \in [0,T) \\
					0 &\geq  \vectorstyle{H}(\vectorstyle{x}(t),\tilde{\vectorstyle{u}}(t),\vectorstyle{z}(t)), ~  \forall t \in [0,T) \\
					0&= \vectorstyle{x}(0) -\vectorstyle{x}(T) \\
					0&= \Psi(\vectorstyle{x}(0))	\end{aligned} \\
            \end{aligned}
	\end{equation}
    where $\vectorstyle{\Phi} \in [ \underline{\vectorstyle{\Phi}}, \overline{\vectorstyle{\Phi}} ]$  are bounded decision variables (homotopy parameters) that can be used to gradually increase the complexity of the problem, and $\vectorstyle{s}\in \mathbb{R}_{+}^{n_{\Phi}}$ are positive penalty variables. The differential and algebraic variables are equivalent to those in problem (\ref{eq:ocp}), but we note that the input variables are updated as
    \begin{subequations}
        \begin{align}
            \label{eq:ctrl_inputs_homotopy}
		  \tilde{\vectorstyle{u}} &= \begin{pmatrix}
		      \vectorstyle{u}^\top &  \hat{\vectorstyle{u}}^\top
		  \end{pmatrix}^\top \\
    \hat{\vectorstyle{u}} &= 
    \begin{pmatrix} 
            \vectorstyle{f}_{\text{fict}}^{\mathcal{B},\top}&
            \vectorstyle{\tau}_{\text{fict}}^{\mathcal{B},\top}&
            {f}_{\text{aircraft,prop}}^\mathcal{B}
            \end{pmatrix}^\top
        \end{align}
    \end{subequations}
    Here, $\vectorstyle{f}_{\text{fict}}^{\mathcal{B}}$ and $\vectorstyle{\tau}_{\text{fict}}^{\mathcal{B}}$ are fictitious forces and torques that are used in early homotopy stages and that are assumed to act directly on the aircraft body instead of the aerodynamic force and torques. A non-physical propulsion/brake force $\vectorstyle{f}_{\text{aircraft,prop}}^\mathcal{B}$ is also applied in $x$-axis of body reference frame. Accordingly, the optimization parameters in (\ref{eq:ocp_homotopy}) are defined as $\tilde{\vectorstyle{w}} = \begin{pmatrix}
        \vectorstyle{x}^\top & \tilde{\vectorstyle{u}}^\top & \vectorstyle{z}^\top \end{pmatrix}^\top$.
    
    For this study, three homotopy stages are used (\( n_\Phi = 3 \)).  In the first homotopy step, fictitious forces and torques are gradually removed from the initial solver, restoring the original aerodynamics. This transition is facilitated by relying on the propulsion/braking force. In the second stage, the propulsion or braking force is set to zero. The DAE $\vectorstyle{F}_{\Phi}$ in (\ref{eq:ocp_homotopy}) becomes similar to $\vectorstyle{F}$, except that the aerodynamic forces applied to the aircraft in (\ref{eq:aircraft}) are replaced by 
    \begin{equation}
        \label{eq:homotopy_stage_1}
        \begin{multlined}
                \begin{pmatrix} {\vectorstyle{f}}_{\text{aircraft},\text{aero}}^\mathcal{B} \\ {\vectorstyle{\tau}}_{\text{aircraft},\text{aero}}^\mathcal{B} 
        \end{pmatrix} \leftarrow (1 - \Phi_1)  \begin{pmatrix} {\vectorstyle{f}}_{\text{aircraft},\text{aero}}^\mathcal{B} \\ {\vectorstyle{\tau}}_{\text{aircraft},\text{aero}}^\mathcal{B}
        \end{pmatrix} \\+ \Phi_1 \begin{pmatrix} \vectorstyle{f}_{\text{fict}}^\mathcal{B} \\ \vectorstyle{\tau}_{\text{fict}}^\mathcal{B} 
        \end{pmatrix}
         + \Phi_2 \; {f}_{\text{aircraft,prop}}^\mathcal{B} \vectorstyle{E}_1^\mathcal{B}
        \end{multlined}
    \end{equation}

    To expedite solving the OCP with the most accurate model, the cost function is initially designed to minimize aircraft angular acceleration and side slip, keeping the aircraft close to its initial path by the end of the second homotopy stage. In the final homotopy stage, the cost function shifts its focus to the power harvesting mode. The cost function $c_{\phi}$ is defined as 
    \begin{equation}
                \label{eq:cost_homotopy}
        \begin{multlined}
            c_{\phi}(\vectorstyle{x}(t), \vectorstyle{u}(t), \vectorstyle{z}(t),\vectorstyle{\phi}) =  \\ -(1-\Phi_3) P + \|\vectorstyle{\xi}\|_\Gamma^2  + \Phi_3\| \vectorstyle{x}(t) - \vectorstyle{x}_0(t) \|_\Lambda^2
        \end{multlined}
    \end{equation}
    where $\Gamma$ and $\Lambda$ are constant diagonal weight matrices, and $\vectorstyle{x}_0(t)$ is the initial state trajectory calculated through the initial path generation.
    
    Now, the solution to the original optimization problem is obtained by solving \(2n_\Phi \) consecutive NLPs, with each solution serving as the warm start for the subsequent problem. In this work, we use IPOPT \cite{biegler2009large}. As IPOPT solves a relaxed version of the Karush-Kuhn-Tucker (KKT) conditions, warm starting is less efficient. Consequently, a Penalty-based Interior-Point Homotopy (PIPH) method \cite{RN5}, which set appropriate values for the internal IPOPT hyper-parameters according to the current requirements of the homotopy step, is used.
    
    \section{Results and Discussion}
    \label{sec_4}
    The simulation results and optimal trajectories derived from the framework introduced in section \ref{sec_3} using the model of section \ref{sec_2} are presented here. Our focus lies on the megawatt-scale AWES with a two-fuselage aircraft having a wing span of $42.5m$ meters, as outlined in \cite{RN17}. We investigated the lemniscate and circular paths for a logarithmic wind shear profile with speed of $15 m/s$. To understand how a flexible tether affects the dynamics of AWES from an optimization point of view, we simplify the quasi-static approach by assuming the tether is rigid.
 
\begin{figure}[ht!]
                \centering
                \includegraphics[trim={3cm 8.5cm 3cm 9cm}, clip, width=.5\textwidth]{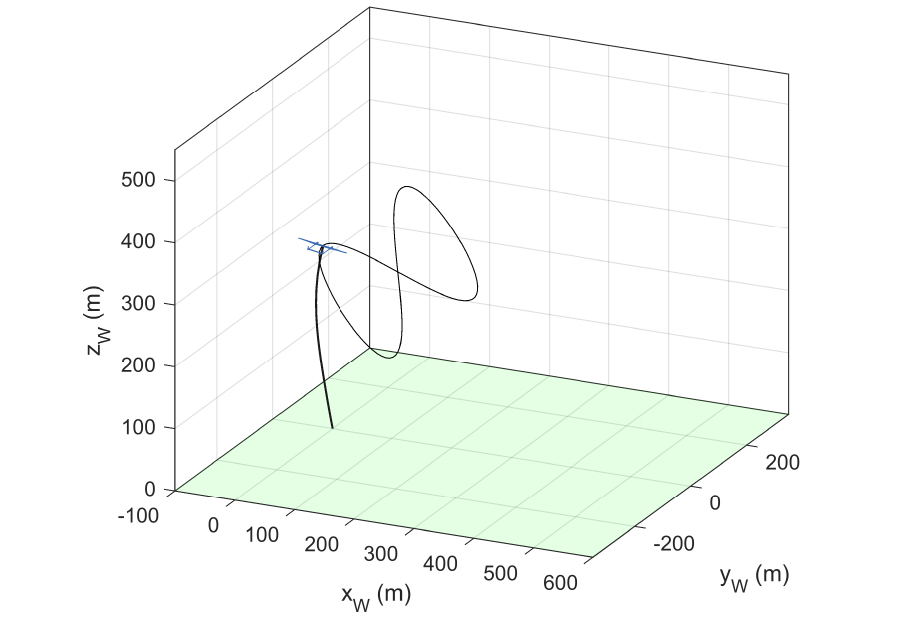}
                \caption{Visualization of the initial lemniscate path.} 
                \label{fig:lemniscate_init_15_3_loops}
            \end{figure}

    The following results and optimal trajectories are derived from the optimization problem (\ref{eq:ocp_homotopy}), where we provide well-defined initial circular and lemniscate trajectories. The flight time is estimated by dividing the approximate total circumference of the initial flight path by a user-defined flight speed parameter. For example, the Fig.~\ref{fig:lemniscate_init_15_3_loops} only shows a three-pumping-cycle lemniscate trajectory, in which the aircraft flies a lemniscate three times within $90s$. We discretize OCP into $N_{t} = 270$ shooting nodes such that time step $h = \frac{T_{f}}{N_{t}} = 0.3 s$. Since determining an accurate initial flight trajectory for the OCP is a non-trivial task, a simple flight with a constant tether length is considered. Repeating this maneuver serves as a suitable candidate for generating trajectories with multiple loops. 
                 
    Considering the AWES model with both flexible and rigid tethers is inherently index-1, $\Psi(\vectorstyle{x}(0))$ in (\ref{eq:ocp}) is considered as: 
    \begin{equation}
        \begin{aligned}
            \label{eq:wind}
            \Psi(\vectorstyle{x}(0)) = \begin{pmatrix}
            \vectorstyle{g}(\vectorstyle{x}_{0},\vectorstyle{z}_{0}) \\   
            \vectorstyle{p}_\text{aircraft}^\mathcal{W}(t=0) - \vectorstyle{p}_{0}
            \end{pmatrix}
        \end{aligned}
    \end{equation}
   where the first row is a cold start \cite{RN9} for the solution, and the second one can be used to constrain the initial point of the optimal path to a specific location, which could be either the starting point of the initial path or a point where the take-off maneuver is supposed to end. Constraining the initial position of the aircraft or any other quantities at this point is optional and can be implemented if there are specific physical or operational requirements. The critical assumptions of non-singular $\frac{\partial{\vectorstyle{g}}}{\partial{\vectorstyle{z}}}$ are achieved if the tether remains under tension $z_3>0$, and none of the tether lumped masses collide with the ground $z_2>0$.
     
    To investigate the effect of the flexible tether model on path planning results, we compare our findings with those obtained using a rigid tether model. Previous studies on the optimization and path planning of rigid wing AWES \cite{RN15,RN47,RN5,RN13} utilized a non-minimal coordinate model and considered the tether force as a constraint. However, it is impossible to incorporate that constraint within our model. Therefore, we introduced the following model, which allows us to contrast rigid and flexible tether modeling in path planning and optimal control problems.

    \subsection{Rigid tether model (simplified quasi-static approach)}
    
    To model a rigid tether, we simplified the flexible tether by neglecting the drag and gravitational forces acting on the tether's lumped masses. By eliminating these forces, which cause tether sag, we effectively made the tether rigid. Consequently, the force at the winch directly points to the aircraft's position. Therefore, equation (\ref{eq:tether_forces}) simplifies to 
    \begin{equation}
            \label{eq:tether_forces_ablation}
		m_j \ddot{\vectorstyle{p}}_j = \vectorstyle{f}_{j-1} - \vectorstyle{f}_{j}
	\end{equation}
    
    The roots of equation (\ref{eq:tether_algebraic}) represent the tether force magnitude and direction at the winch, which can be numerically found using the direct multiple shooting method, similar to the flexible tether case.
    To estimate the tether's drag force act on the aircraft, we make use of a method similar to \cite{RN8, argatov2013efficiency}, in which the drag force for the whole tether length is given by 
    \begin{equation}
            \label{eq:tether_drag_ablation}
		\vectorstyle{f}_{\text{tether}}^\text{drag} = -\frac{1}{8} \rho
            C_D l_\text{tether} d_\text{tether}\|\vectorstyle{v}_{a,\perp}\| \vectorstyle{v}_{a,\perp}
	\end{equation}
    where $l_\text{tether}$ is the whole tether length. The drag and gravitational forces are then applied to the aircraft by adding them to (\ref{eq:forces}), in order to have a fair comparison between flexible and rigid tether models.

    \subsection{Simulation results}
\begin{figure}[hbt!]
                \centering
                \includegraphics[trim={3.5cm 8.5cm 3.5cm 8.5cm}, clip,  width=.5\textwidth]{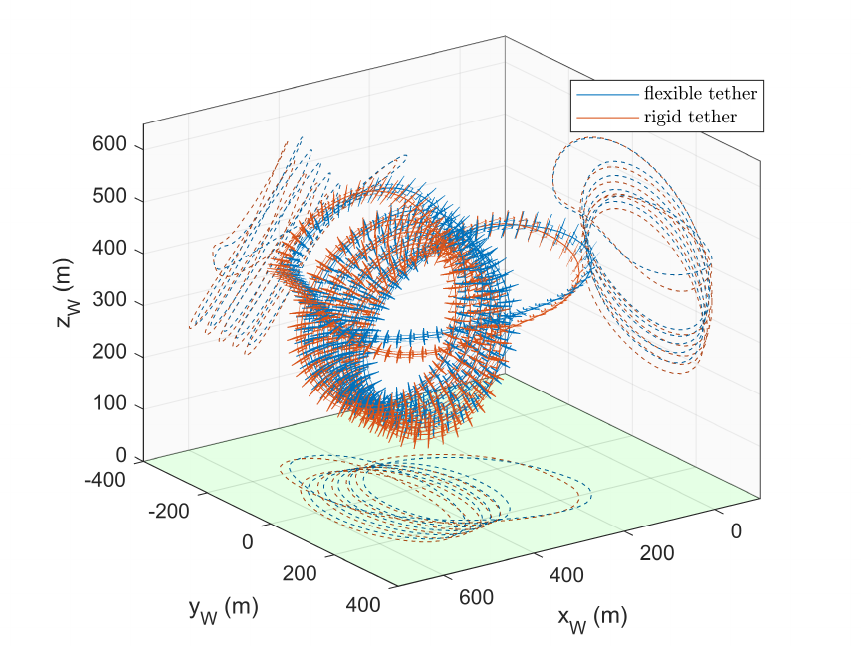}
                \caption{Visualization of optimal trajectory with six circular pumping cycles for quasi-static and simplified  tether models.} 
                \label{fig:circle_normal_ablation}

        \end{figure}

\begin{figure}[hbt!]      
                \centering
                \includegraphics[trim={3.5cm 8.5cm 3.5cm 8.5cm}, clip,  width=.5\textwidth]{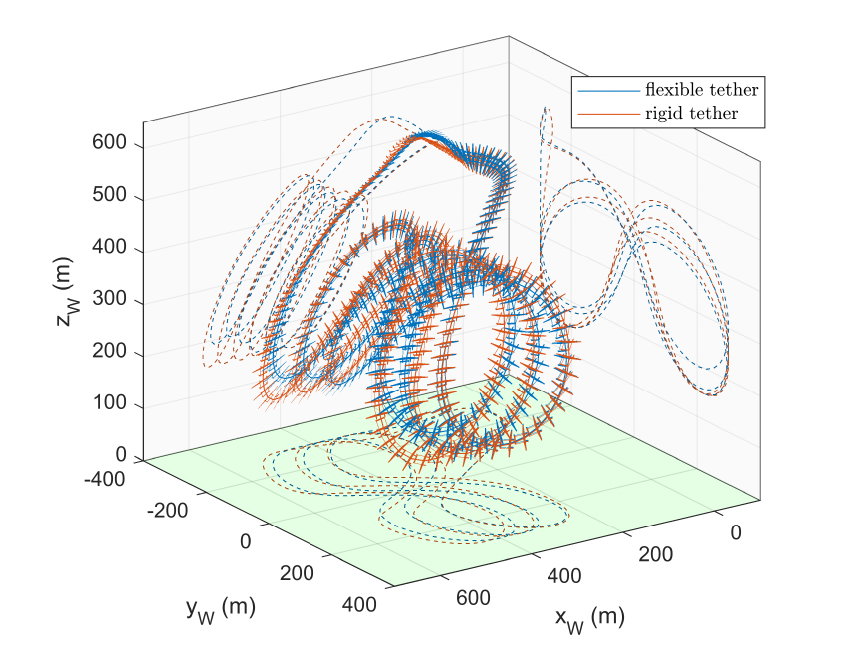}
                \caption{Visualization of optimal trajectory with three lemniscate pumping cycles for quasi-static and simplified  tether models.} 
                \label{fig:lemniscate_normal_ablation}
        \end{figure}

    In this section, we discuss the generated optimal trajectories for both flexible and rigid tethers. As it is depicted in Fig.~\ref{fig:circle_normal_ablation} and Fig.~\ref{fig:lemniscate_normal_ablation}, the optimal paths for rigid and flexible tethers look similar. However, some differences are observed both in the traction and retraction phases. These discrepancies are more pronounced during the retraction and the transition between the traction and retraction phases.

    Despite the resemblance between the optimal trajectories for models with rigid and flexible tethers, there are notable differences in the control input signals and the estimated harvested power. In other words, the optimizer aims to maximize power by producing maneuvers that increase lift. Therefore, while it is expected that the flying trajectories for both rigid and flexible tethers remain similar, the differences in tether dynamics and their effects on the aircraft cause changes in the control inputs.

\begin{figure}[hbt!]
                    \centering
                    \includegraphics[trim={3cm 8.5cm 3cm 8.5cm}, clip,  width=.5\textwidth]{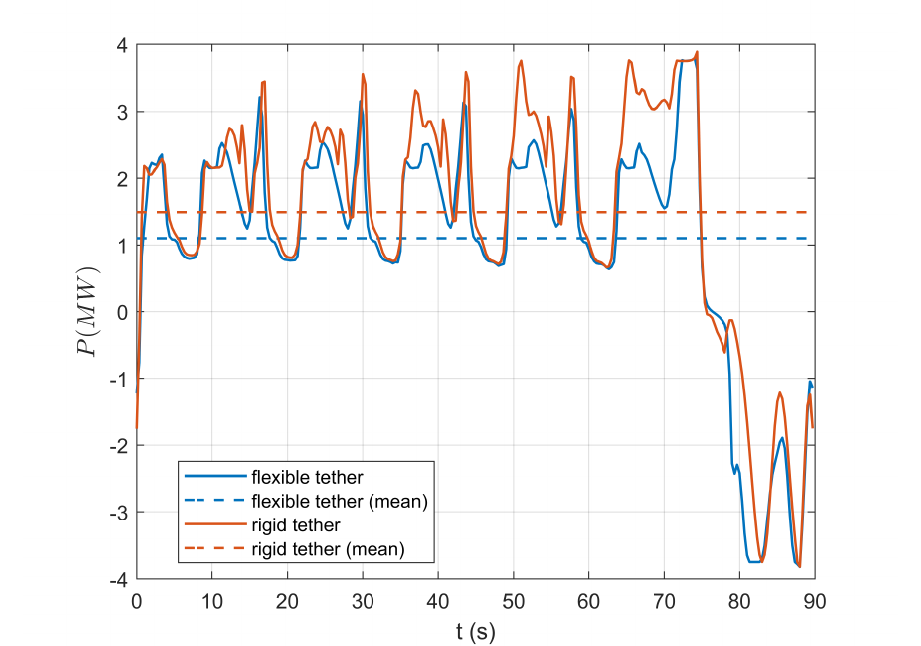}
                    \caption{The mechanical power and its corresponding average power for
         optimal circular trajectories with quasi-static and simplified tether models.} 
                    \label{fig:Power_0}
         \end{figure}

\begin{figure}[hbt!]
                    \centering
                   \includegraphics[trim={3cm 8.5cm 3cm 8.5cm}, clip,  width=.5\textwidth]{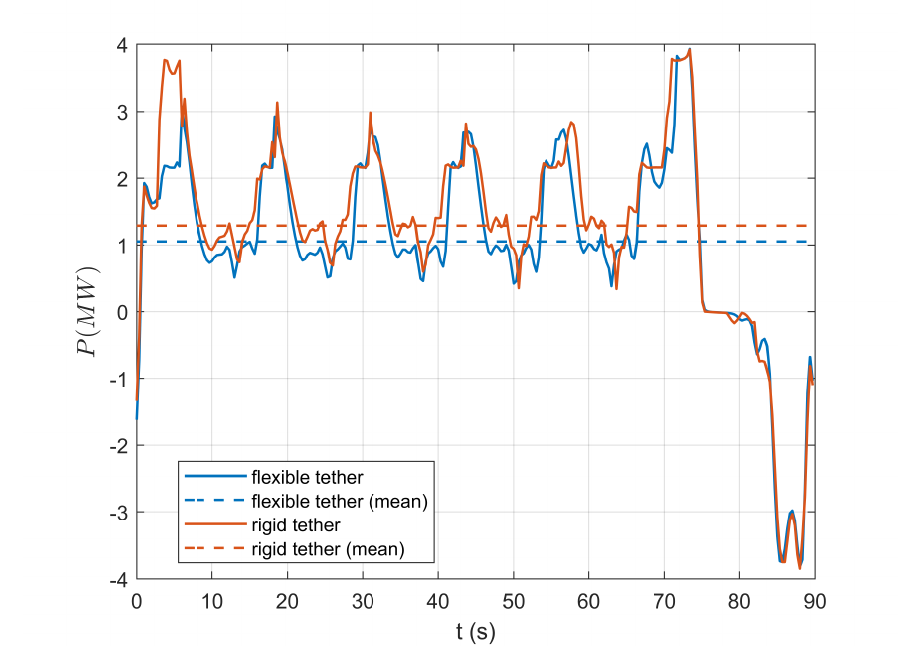}
                   \caption{The mechanical power and its corresponding average power for
        optimal lemniscate trajectories with quasi-static and simplified tether models.} 
                    \label{fig:Power_8}
        \end{figure}

    Although \cite{RN8} claimed that the rigid tether model can accurately estimate power for one pumping cycle, the results shown in Fig.~\ref{fig:Power_0} and Fig.~\ref{fig:Power_8} demonstrate that it underestimates the effect of tether sag on the harvested power for both circular and lemniscate paths, in which the power in the pumping cycles of the model with a rigid tether is higher than that of the model with a flexible tether.
    
        \begin{figure}[hbt!]
                \centering
                \includegraphics[trim={3cm 8.5cm 3cm 8.5cm}, clip,  width=.5\textwidth]{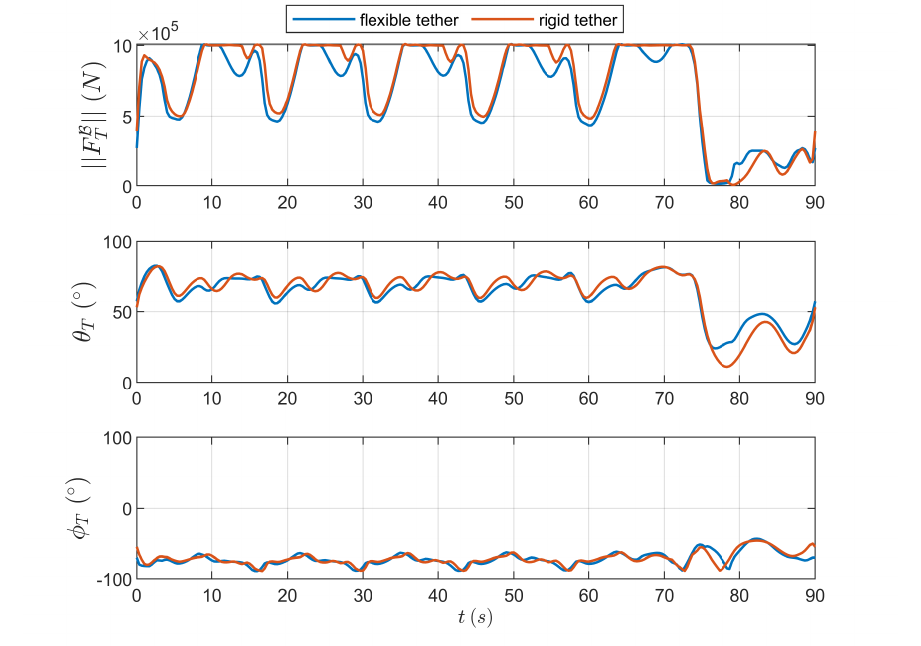}
                \caption{Tether force characteristics in body frame for the circular path.} 
                \label{fig:tether_force_0}
        \end{figure}
        
        \begin{figure}[hbt!]
                \centering
               \includegraphics[trim={3cm 8.5cm 3cm 8.5cm}, clip,  width=.5\textwidth]{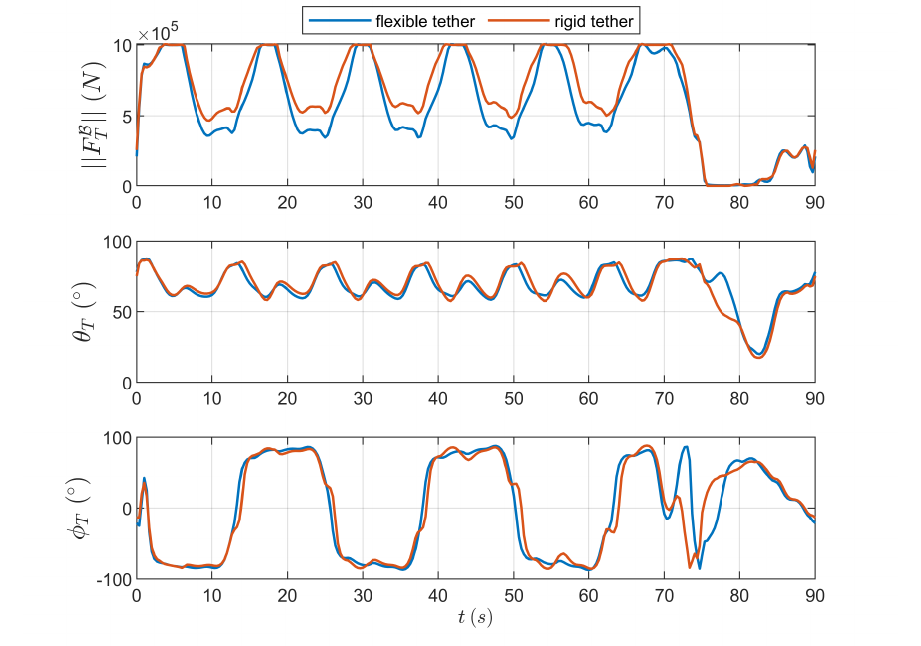}
                \caption{Tether force characteristics in body frame for the lemniscate path.} 
            \label{fig:tether_force_8}
        \end{figure}
    
    The tether force acting on the aircraft can be expressed based on its magnitude and two angles, which determine its direction in the body reference frame (Fig.~\ref{fig:tether_collision}). We previously described how we used one of those angles to avoid tether collision with aircraft. In Fig.~\ref{fig:tether_force_0} and Fig.~\ref{fig:tether_force_8}, we quantify how tether sag might affect the direction of the tether force acting on the aircraft. The difference in angles during the retraction phase, for times after $75s$, represents the tether sag captured in the flexible tether model.      
    
\begin{figure}[hbt!]
                \centering
                \includegraphics[trim={4.0cm 8.5cm 4.5cm 8.5cm}, clip,  width=.5\textwidth]{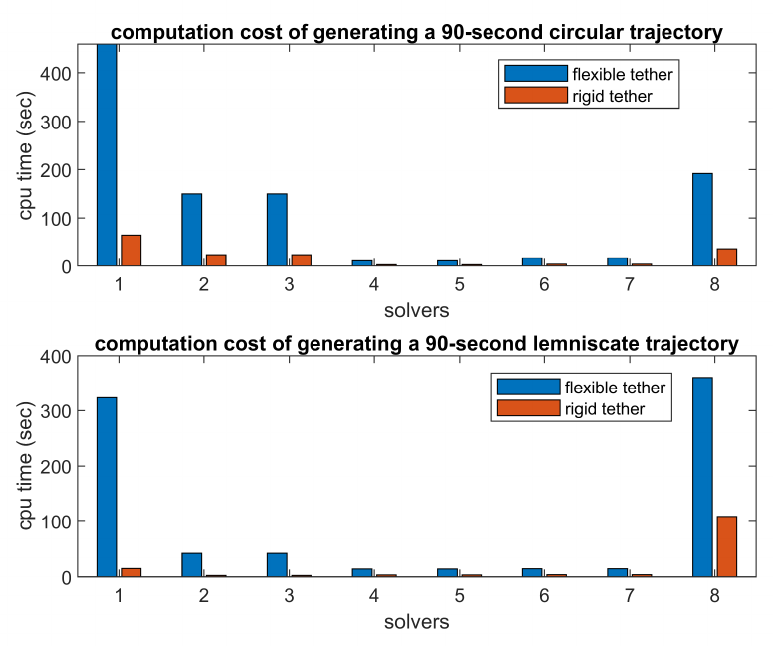}
                \caption{Computation cost of 90-second trajectories with flexible and rigid tether models.} 
                \label{fig:computational_cost_08}
            \end{figure}

    Addressing computational effort is essential for optimization purposes. The CPU time required to solve the optimal control problem on a standard laptop with a 1.8 GHz Intel CPU and 16 GB RAM for circular and lemniscate trajectories is presented below. We assumed ten segments on the tether, which is sufficient to show the tether sag while keeping the computation time reasonable. A higher number of segments increases the model's non-linearity and complexity, causing the computation time for the flexible tether to be higher than that for the rigid one.
    As can be seen in Fig.~\ref{fig:computational_cost_08}, the computation time for finding an optimal lemniscate trajectory is much higher than for the circular one. The more complicated and less well-defined forces and torques in lemniscate trajectories may be the reason for the increased computational effort.

     Solver numbers one to eight are the homotopy solvers used for three homotopy stages introduced earlier. Solver numbers one and eight are the initial and final solvers, respectively. The intermediate solvers are divided into three homotopy stages: the first stage includes solvers two and three, which remove fictitious forces and torques; the second stage includes solvers four and five, which eliminate the propulsion force; and the third stage includes solvers six and seven, during which the cost function is adjusted for power production.

\section{Conclusion}
    \label{sec_5}
    This work presents an optimal path planning method to evaluate the performance of an \textsl{airborne wind energy system} (AWES), which includes a flexible tether, a winch, and a rigid-wing aircraft. By the proposed formulation, we generate paths with multiple pumping cycles using a homotopy approach. Unlike other research studies that employ an index-3 \textsl{differential-algebraic system of equations} (DAE) form and solve it using index reduction techniques, we employ a semi-explicit index-1 DAE form. This allows us to add complementary constraints at the initial point to enhance the trajectory characteristics. Moreover, the quasi-static approach permits the use of a flexible tether model within an optimization problem without significantly increasing the number of states, as might occur with other lumped mass models used in the literature.
    
    To better understand how the trajectories generated with a flexible tether model differ from the ones generated with a rigid tether, a simplified version of the quasi-static approach was introduced to imitate the behavior of a rigid tether. Since the tether's drag force was dropped as part of the simplification, we made use of a tether drag force estimation. The simulation results reveal that, although the optimal trajectories for an AWES with rigid and flexible tethers appear similar, the direction of the force exerted by the tether on the aircraft can differ significantly during the retraction phase. This difference leads to more pronounced mismatches between the optimal trajectories of the rigid and flexible tethers during the retraction phase compared to the pumping cycles. This can cause an overestimation of harvested power for one pumping cycle by up to 33\%.

\section*{Funding Sources}
 With the support of Energy Transition Fund of FPS Economy, project BORNE (\textsl{Belgian Offshore aiRborne wiNd Energy}).

\bibliography{sample}

\end{document}